\documentclass{article}

\oddsidemargin=\evensidemargin
\catcode`\@=11
\long\def\@makefntext#1{
\protect\noindent \hbox to 3.2pt {\hskip-.9pt  
$^{{\eightrm\@thefnmark}}$\hfil}#1\hfill}		

\def\ps@myheadings{\let\@mkboth\@gobbletwo		
\def\@oddhead{\hbox{}
\rightmark\hfil\eightrm\thepage}   
\def\@oddfoot{}\def\@evenhead{\eightrm\thepage\hfil
\leftmark\hbox{}}\def\@evenfoot{}
\def\sectionmark##1{}\def\subsectionmark##1{}}

\catcode`@=11		      		     
\def\ps@plain{\let\@mkboth\@gobbletwo
     \def\@oddhead{}\def\@oddfoot{\eightrm\hfil\thepage
     \hfil}\def\@evenhead{}\let\@evenfoot\@oddfoot}

\renewcommand{\thefootnote}{\fnsymbol{footnote}}

\newcounter{sectionc}\newcounter{subsectionc}\newcounter{subsubsectionc}
\renewcommand{\section}[1] {\vspace{12pt}\addtocounter{sectionc}{1} 
\setcounter{subsectionc}{0}\setcounter{subsubsectionc}{0}\noindent 
	{\tenbf\thesectionc. #1}\par\vspace{5pt}}
\renewcommand{\subsection}[1] {\vspace{12pt}\addtocounter{subsectionc}{1} 
	\setcounter{subsubsectionc}{0}\noindent 
	{\bf\thesectionc.\thesubsectionc. 
	{\kern1pt \bfit #1}}\par\vspace{5pt}}
\renewcommand{\subsubsection}[1] {\vspace{12pt}
	\addtocounter{subsubsectionc}{1}
	\noindent
	{\tenrm\thesectionc.\thesubsectionc.\thesubsubsectionc.	{\kern1pt 
	\it #1}}\par\vspace{5pt}}
\newcommand{\nonumsection}[1] {\vspace{12pt}\noindent{\tenbf #1}
	\par\vspace{5pt}}

\newcounter{appendixc}
\newcounter{subappendixc}[appendixc]
\newcounter{subsubappendixc}[subappendixc]

\renewcommand{\appendix}[1] {\vspace{12pt}	
	\refstepcounter{appendixc}		
	\setcounter{figure}{0}
	\setcounter{table}{0}
	\setcounter{lemma}{0}
	\setcounter{theorem}{0}
	\setcounter{corollary}{0}
	\setcounter{definition}{0}
	\setcounter{equation}{0}
	\renewcommand{\thefigure}{\Alph{appendixc}.\arabic{figure}}
	\renewcommand{\thetable}{\Alph{appendixc}.\arabic{table}}
	\renewcommand{\theappendixc}{\Alph{appendixc}}
	\renewcommand{\thelemma}{\Alph{appendixc}.\arabic{lemma}}
	\renewcommand{\thetheorem}{\Alph{appendixc}.\arabic{theorem}}
	\renewcommand{\thedefinition}{\Alph{appendixc}.\arabic{definition}}
	\renewcommand{\thecorollary}{\Alph{appendixc}.\arabic{corollary}}
	\renewcommand{\theequation}{\Alph{appendixc}.\arabic{equation}}
	\noindent{\tenbf Appendix \theappendixc #1}\par\vspace{5pt}}

\newcommand{\textlineskip}{\baselineskip=13pt}
\newcommand{\smalllineskip}{\baselineskip=10pt}

\newcommand{\copyrightheading}[1]
	{\vspace*{-2.5cm}\smalllineskip{\flushleft
	{\footnotesize Journal of Knot Theory and Its Ramifications #1}\\
   	{\footnotesize \copyright\kern2pt World Scientific 
         Publishing Company}\\
         }}

\def\abstracts#1#2#3#4{{
	\centering{\begin{minipage}{4.5in}\footnotesize\baselineskip=10pt
	\centerline{ABSTRACT} 
	\parindent=15pt #1\par 
	\parindent=15pt #2\par
	\parindent=15pt #3\par
	\parindent=15pt #4\par
	\end{minipage}}\par}} 

\renewenvironment{thebibliography}[1]
	{\frenchspacing
	 \ninerm\baselineskip=11pt
	 \begin{list}{[\arabic{enumi}]}
	{\usecounter{enumi}\setlength{\parsep}{0pt}
	 \setlength{\leftmargin 13.7pt}{\rightmargin 0pt} 
	 \setlength{\itemsep}{0pt} \settowidth
	{\labelwidth}{[#1]}\sloppy}}{\end{list}}

\newcounter{itemlistc}
\newcounter{romanlistc}
\newcounter{alphlistc}
\newcounter{arabiclistc}

\newcommand{\fcaption}[1]{
        \refstepcounter{figure}
        \setbox\@tempboxa = \hbox{\footnotesize Fig.~\thefigure. #1}
        \ifdim \wd\@tempboxa > 5in
           {\begin{center}
        \parbox{5in}{\footnotesize\smalllineskip Fig.~\thefigure. #1}
            \end{center}}
        \else
             {\begin{center}
             {\footnotesize Fig.~\thefigure. #1}
              \end{center}}
        \fi}

\newcommand{\tcaption}[1]{
        \refstepcounter{table}
        \setbox\@tempboxa = \hbox{\footnotesize Table~\thetable. #1}
        \ifdim \wd\@tempboxa > 5in
           {\begin{center}
        \parbox{5in}{\footnotesize\smalllineskip Table~\thetable. #1}
            \end{center}}
        \else
             {\begin{center}
             {\footnotesize Table~\thetable. #1}
              \end{center}}
        \fi}

\def\pmb#1{\setbox0=\hbox{#1}
	\kern-.025em\copy0\kern-\wd0
	\kern.05em\copy0\kern-\wd0
	\kern-.025em\raise.0433em\box0}

\def\fnt#1#2{\footnotetext{\kern-.3em
	{$^{\mbox{\scriptsize #1}}$}{#2}}}

\def\fpage#1{\begingroup
\voffset=.3in
\thispagestyle{empty}\begin{table}[b]\centerline{\footnotesize #1}
	\end{table}\endgroup}

\def\runninghead#1#2{\pagestyle{myheadings}
\markboth{{\protect\footnotesize\it{\quad #1}}\hfill}
{\hfill{\protect\footnotesize\it{#2\quad}}}}

\font\tenrm=cmr10
 
\font\tenbf=cmbx10
\font\bfit=cmbxti10 at 10pt
\font\ninerm=cmr9

\font\eightrm=cmr8

\newtheorem{theorem}{Theorem}   

\newtheorem{lemma}{Lemma}

\newtheorem{corollary}{Corollary}
\def\@begintheorem#1#2{\trivlist	
	\item[\hskip\labelsep{\bf #1\ #2.}]} 
\def\@opargbegintheorem#1#2#3{\trivlist
	\item[\hskip\labelsep{\bf #1\ #2\ (#3).}]}

    	{\setcounter{itemlistc}{0}		
	 \begin{list}{$\bullet$}		
	{\usecounter{itemlistc}			
	 \leftmargin10pt	       
	 \setlength{\parsep}{0pt}
	 \setlength{\itemsep}{0pt}     
	}}{\end{list}}

	{\setcounter{romanlistc}{0}		
	 \begin{list}{$($\roman{romanlistc}$)$}	
	{\usecounter{romanlistc}		
	 \leftmargin18pt 
	 \setlength{\parsep}{0pt}
	 \setlength{\itemsep}{0pt}	
	 \settowidth{\labelwidth}{#1}                          
	}}{\end{list}}

	{\setcounter{enumii}{0}			
	 \begin{list}{$($\alph{enumii}$)$}	
	{\usecounter{enumii}			
	 \leftmargin18pt		
	 \setlength{\parsep}{0pt}
	 \setlength{\itemsep}{0pt}	
	 \settowidth{\labelwidth}{#1}                          
	}}{\end{list}}

\def\qed{\hbox{${\vcenter{\vbox{			
   \hrule height 0.4pt\hbox{\vrule width 0.4pt height 6pt
   \kern5pt\vrule width 0.4pt}\hrule height 0.4pt}}}$}}

\renewcommand{\thefootnote}{\fnsymbol{footnote}}  

\def\theequation{\thesectionc.\arabic{equation}}  
\begin{document}
\setlength{\textheight}{7.7truein}  

\runninghead{Impossibility of obtaining split links}
{Impossibility of obtaining split links}

\normalsize\textlineskip
\thispagestyle{empty}
\setcounter{page}{1}

\copyrightheading{}		    

\vspace*{0.88truein}

\fpage{1}
\centerline{\bf IMPOSSIBILITY OF OBTAINING SPLIT LINKS}
\baselineskip=13pt
\centerline{\bf FROM SPLIT LINKS VIA TWISTINGS
}
\vspace*{0.37truein}
\centerline{\footnotesize MAKOTO OZAWA
\footnote{The author was supported in part by Fellowship of
the Japan Society for the Promotion of Science for Japanese Junior Scientists.}
}
\baselineskip=12pt
\centerline{\footnotesize\it Department of Mathematics, School of Education,}
\centerline{\footnotesize\it Waseda University, Nishiwaseda 1-6-1, Shinjuku-ku,}
\centerline{\footnotesize\it Tokyo 169-8050, Japan}
\baselineskip=10pt
\centerline{\footnotesize\it ozawa@musubime.com}


\vspace*{0.225truein}
\vspace*{1truein}

\vspace*{0.21truein} 
\abstracts{We show that if a split link is obtained from a split link $L$ in $S^3$ by $1/n$-Dehn surgery along a trivial knot $C$, then the link $L\cup C$ is splittable. That is to say, it is impossible to obtain a split link from a split link via a non-trivial twisting. As its corollary, we completely determine when a trivial link is obtained from a trivial link via a twisting.
}{}{}{}



\vspace*{1pt}\textlineskip	

\section{Introduction}	
\vspace*{-0.5pt}

Let $L$ be a link in $S^3$ and $C$ a trivial knot in $S^3$ missing $N(L)$. Then we can get a new link $L^*$ in $S^3$ as the image of $L$ after doing $1/n$-Dehn surgery along $C$. We say that $L^*$ is obtained from $L$ by an {\it $n$-twisting} along $C$. In this paper, we consider the following problem.

\vspace*{2mm}
\noindent{\bf Problem.}
Is it possible that both $L$ and $L^*$ are splittable?
\vspace*{2mm}

For this problem, it is reasonable to make the following definition.
An $n$-twisting is said to be {\it non-trivial} if $n\ne 0$ and the link $L\cup C$ is non-splittable.
Then our result is stated as follows.

\vspace*{1mm}
\begin{theorem}
It is impossible to obtain a split link from a split link by a non-trivial twisting.
\end{theorem}
\vspace*{1mm}

Next, we consider when a trivial link is obtained from a trivial link by an $n$-twisting. For a trivial knot, this problem has been solved as follows.

\vspace*{1mm}
\begin{theorem} $($\cite{M}, \cite{KMS}$)$
Suppose that a trivial knot $K^*$ is obtained from a trivial knot $K$ by an  $n$-twisting along $C$. Then one of the following conclusions holds.
\par
(1) The link $K\cup C$ is a trivial link.
\par
(2) The link $K\cup C$ is a Hopf link.
\par
(3) The link $K\cup C$ is a torus link of type $(4,-2)$ or $(4,2)$, and $n=1$ or $-1$ respectively.
\end{theorem}
\vspace*{1mm}

By Theorems 1 and 2, we obtain the next corollary.

\vspace*{1mm}
\begin{corollary}
Suppose that a trivial link $L^*=K_1^*\cup \ldots \cup K_l^*$ is obtained from a trivial link $L=K_1\cup \ldots \cup K_l$ by an $n$-twisting along $C$. Then one of the following conclusions holds.
\par 
(1) The link $L\cup C$ is a trivial link.
\par 
(2) The link $L\cup C$ is a split union of a Hopf link $K_i\cup C$ and a trivial link $L-K_i$ for some $i\in \{1, \ldots , l\}$ 
\par 
(3) The link $L\cup C$ is a split union of a torus link $K_i\cup C$ of type $(4,-2)$ or $(4,2)$ and a trivial link $L-K_i$ for some $i\in \{1,\ldots ,l\}$, and $n=1$ or $-1$ respectively.
\end{corollary} 
\vspace*{1mm}

\section{Preliminaries}

In this section, we prepare some lemmas for Theorem 1.
All manifolds are assumed to be compact and orientable, and any srufaces in a 3-manifold are assumed to be properly embedded and in general position.

Let $M$ be a 3-manifold, and $F_1$ and $F_2$ two surfaces in $M$. Let $\hat{F_1}$, $\hat{F_2}$ be the closed surfaces obtained by capping off $\partial F_1$, $\partial F_2$ with disks. Then, for $\alpha \in \{ 1,2\}$, one defines a graph $G_{\alpha}$ in $\hat{F_{\alpha}}$, where the edges of $G_{\alpha}$ correspond to the arc components of $F_1\cap F_2$, and the vertices to the components of $\partial F_{\alpha}$. Recall that a {\it 1-sided face} in a graph is a disk face with exactry one edge in its boundary.

Recall that if $M$ is a 3-mainfold with torus boundary and $\gamma$ is a slope on $\partial M$, then $M(\gamma )$ denotes the closed maifold obtained by attaching a solid torus $J$ to $M$ so that the boundary of a meridian disk of $J$ has slope $\gamma$ on $\partial M$. Recall also that if $\gamma _1,\gamma _2$ are two slopes on $\partial M$, then $\Delta (\gamma _1,\gamma _2)$ denotes the minimal geometric intersection number of $\gamma _1$ and $\gamma _2$.

The following lemma will be needed for Theorem 1.

\vspace*{1mm}
\begin{lemma}
Let $M$ be a 3-manifold with torus boundary and let $F_1,F_2$ be planar surfaces in $M$ with boundary slopes $\gamma _1,\gamma _2$. Suppose that the graphs $G_1,G_2$ contain no 1-sided faces, and that $\Delta (\gamma _1,\gamma_2)\ge 1$. Then either the first homology groups $H_1(M(\gamma _1))$ or $H_1(M(\gamma_2))$ has a torsion.
\end{lemma} 
\vspace*{1mm}

\noindent{\it Proof.}
If $\Delta (\gamma _1,\gamma _2)\ge 2$, then Lemma 1 follows \cite[lemma 2.2]{GL1}.
Otherwise, by  \cite[Proposition 2.0.1]{GL2}, $G_1$ contains a Scharlemann cycle or $G_2$ represents all $\{ 1,\ldots ,|\partial P_1|\}$-types. In the formar case, $M(\gamma _2)$ has a lens space as a connected summand. In the latter case, by  \cite[Theorem]{P}, $H_1(M(\gamma _1))$ has a torsion. This completes the proof of Lemma 2.
\qed\kern0.8pt

\section{Proof of Theorem 1} 

Suppose that a split link $L^*$ is obtained from a split link $L$ by a non-trivial $n$-twisting along $C$. Let $S$ and $S^*$ be the splitting spheres for $L$ and $L^*$ respectively.
Put $M=S^3-intN(C)$.
We may assume that $C$ intersects $S$ and $S^*$ transversely in the 3-spheres $M(1/0)$ and $M(1/n)$ respectively, and assume that $|C\cap S|$ and $|C\cap S^*|$ are minimal among all 2-spheres isotopic to $S$ and $S^*$ respectively. 
Then, since $L\cup C$ and $L^*\cup C$ are non-splittable, $|C\cap S|$ and $|C\cap S^*|$ are not equal to zero. Put $P_1=S-intN(C)$ and $P_2=S^*-intN(C)$. Then by the minimality of $|C\cap S|$ and $|C\cap S^*|$ and by the irreducibility of $M-L$, $P_1$ and $P_2$ satisfy the hypothesis of Lemma 2. Hence $H_1(M(1/0))$ or $H_1(M(1/n))$ has a torsion, this is impossible.
\qed\kern0.8pt

\medskip

\nonumsection{Acknowledgement}
The author would like to thank Prof. Chuichiro Hayashi for his helpful comments.

\end{document}